\documentclass{CustomArticle}

\usepackage{algorithm}
\usepackage[noend]{algpseudocode}
\usepackage[
  style = trad-abbrv,
  citestyle = numeric-comp,
  sorting = nyt,
  backend = biber
]{biblatex}
\usepackage{booktabs}
\usepackage{threeparttable}
\usepackage{caption}
\usepackage{graphicx}
\usepackage{enumitem}
\usepackage{import}
\usepackage{multirow}
\usepackage{subcaption}
\usepackage{hyperref}
\usepackage[capitalize]{cleveref}

\usepackage{CommonMath}
\usepackage{ProjectHelperMacros}
\usepackage{ProjectStyleTweaks}
\usepackage{ProjectFixes}

\PaperGeneralInfo{
  title={DADA: Dual Averaging with Distance Adaptation},
  date={April 20, 2026}
}
\AddPaperAuthor{
  name={Mohammad Moshtaghifar},
  affiliation={University of British Columbia},
  email={mmoshtaq@cs.ubc.ca}
}
\AddPaperAuthor{
  name={Anton Rodomanov},
  affiliation={CISPA Helmholtz Center for Information Security},
  email={anton.rodomanov@cispa.de}
}
\AddPaperAuthor{
  name={Daniil Vankov},
  affiliation={Arizona State University},
  email={dvankov@asu.edu}
}
\AddPaperAuthor{
  name={Sebastian U. Stich},
  affiliation={CISPA Helmholtz Center for Information Security},
  email={stich@cispa.de}
}
\PaperAbstract{
  We present a novel universal gradient method for solving convex
  optimization problems. Our algorithm---Dual Averaging with Distance Adaptation
  (DADA)---is based on the classical scheme of dual averaging and dynamically
  adjusts its coefficients based on observed gradients and the distance
  between iterates and the starting point, eliminating the need for
  problem-specific parameters. DADA is a universal algorithm that simultaneously
  works for a broad spectrum of problem classes, provided the local growth of
  the objective function around its minimizer can be bounded.
  Particular examples of
  such problem classes are nonsmooth Lipschitz functions, Lipschitz-smooth
  functions, H\"older-smooth functions, functions with high-order Lipschitz
  derivative, quasi-self-concordant functions, and $(L_0,L_1)$-smooth functions.
  Crucially, DADA is applicable to both unconstrained and constrained problems,
  even when the domain is unbounded, without requiring prior knowledge of the
  number of iterations or desired accuracy.
}
\PaperKeywords{
  Convex Optimization, Gradient Methods, Adaptive Algorithms, Dual Averaging,
  Distance Adaption, Universal Methods, Worst-Case Complexity Guarantees
}

\begin{document}

  \PrintTitleAndAbstract

  \begin{abstract}
    We present a novel universal gradient method for solving convex optimization
    problems. Our algorithm---Dual Averaging with Distance Adaptation
    (DADA)---is based on the classical scheme of dual averaging and dynamically
    adjusts its coefficients based on observed gradients and the distance
    between iterates and the starting point, eliminating the need for
    problem-specific parameters. DADA is a universal algorithm whose
    convergence rate adapts to the local behavior of the objective around its
    minimizer, through bounds on its local growth. This leads to a single
    method with explicit, problem-dependent guarantees across a broad range of
    models, including nonsmooth Lipschitz functions, Lipschitz-smooth
    functions, H\"older-smooth functions, functions with high-order Lipschitz
    derivative, quasi-self-concordant functions, and $(L_0,L_1)$-smooth
    functions.
    Crucially, DADA is applicable to both unconstrained and
    constrained problems,
    even when the domain is unbounded, without requiring prior knowledge of the
    number of iterations or desired accuracy.
  \end{abstract}

  \section{Introduction}

Gradient methods are among the most popular and efficient algorithms for
solving optimization problems arising in machine learning, as they are highly
adaptable and scalable across various settings~\cite{Bottou2018OPTMethods}.
Despite their popularity, these methods face a significant challenge of
selecting appropriate hyperparameters,
particularly stepsizes, which are critical to the performance of the algorithm.
Hyperparameter tuning is one of the standard
approaches to address this issue but is a time-consuming and
resource-intensive
process, especially as models become larger and more complex. Consequently,
the cost of training these models has become a significant concern
\cite{sharir2020costtrainingnlpmodels,
  patterson2021carbonemissionslargeneural}.

Typically, line-search techniques have been used to select stepsizes for
optimization methods,
and they are provably efficient for certain function classes, such as
H\"older-smooth problems~\cite{Nesterov2015UniversalGM}. However, in recent
years, several so-called parameter-free algorithms have been developed
which do not utilize line search~\cite{%
  Orabona2017TrainingDN,%
  Cutkosky2018BlackBoxRF,%
  Carmon2022MakingSGDParameterFree,%
  Ivgi2023DoGIS,%
  Khaled2023DoWG,%
  Mishchenko2024Prodigy%
}.
Notably, one strategy involves dynamically adjusting
stepsizes based on estimates of the initial distance to the optimal
solution~\cite{Carmon2022MakingSGDParameterFree,Ivgi2023DoGIS,Khaled2023DoWG}.
Another approach leverages lower bounds on the initial distance combined with
the Dual Averaging (DA)
scheme~\cite{Defazio2023DAdaptation,Mishchenko2024Prodigy}.
However, these methods primarily focus on nonsmooth Lipschitz or, in some
cases, Lipschitz-smooth functions. Some of these methods also come with
additional limitations, such as requiring bounded domain
assumptions~\cite{Khaled2023DoWG} or failing to extend to constrained
optimization problems~\cite{Defazio2023DAdaptation,Mishchenko2024Prodigy}.

To formalize the discussion, we consider the following optimization problem:
\begin{equation}\label{eq:problem}
  f^* \DefinedEqual \min_{x \in Q} f(x),
\end{equation}
where $Q \subseteq \R^d$ is a nonempty closed convex set, and
$\Map{f}{\R^d}{\RealFieldPlusInfty}$ is a proper closed convex function
that is subdifferentiable on $Q$.
We assume that $Q$ is a simple set, meaning that
it is possible to efficiently compute the projection onto $Q$.
We also assume problem~\cref{eq:problem} has a
solution~$x^* \in \Interior \EffectiveDomain f$.
The starting point in our methods is denoted by~$x_0$.

\paragraph{Contributions.}

In this paper, we introduce Dual Averaging with Distance Adaptation (DADA),
a novel universal gradient method for solving~\eqref{eq:problem}.
Building on the classical framework of weighted
DA~\cite{Nesterov2005PrimaldualSM}, DADA incorporates a dynamically adjusted
estimate of $D_0 \DefinedEqual \norm{x_0 - x^*}$, inspired by recent techniques
from~\cite{Ivgi2023DoGIS,Carmon2022MakingSGDParameterFree} and further
developed in~\cite{Khaled2023DoWG}, without requiring prior knowledge of
problem-specific parameters. Furthermore, our approach applies to both
unconstrained problems and those with simple constraints, possibly with
unbounded domains. This makes DADA a powerful tool across a wide range of
applications.

We start, in \cref{sec:preliminaries}, by presenting our
method and outline
its foundational structure based on the DA
scheme~\cite{Nesterov2005PrimaldualSM}. Our main theoretical result,
\cref{thm:main-theorem}, establishes convergence guarantees for a
broad range of function classes.

To demonstrate the versatility and effectiveness of DADA, in
\cref{sec:applications}, we provide complexity estimates across several
interesting function classes:\
nonsmooth Lipschitz functions, Lipschitz-smooth functions, H\"older-smooth
functions, quasi-self-concordant (QSC) functions,
functions with Lipschitz high-order derivative,
and $(L_0, L_1)$-smooth functions.
These results underscore DADA's ability to deliver competitive performance
without knowledge of class-specific parameters.

\paragraph{Related work.}

\begin{table}[tb]
  \centering
  \resizebox{0.95\textwidth}{!}{
    \begin{threeparttable}
      \begin{tabular}{@{}lccccc@{}}
        \toprule
        Method                                                   & Universal                       & Constraints & Unbounded domain & No search &
        Stochastic                                                                                                                                       \\
        \midrule
        DoG~\cite{Ivgi2023DoGIS}                                 & \xmark                          & \cmark      & \cmark           & \cmark    & \cmark \\
        DoWG~\cite{Khaled2023DoWG}                               & \xmark                          & \cmark      & \xmark           & \cmark    & \xmark \\
        Bisection Search~\cite{Carmon2022MakingSGDParameterFree} & \xmark                          & \cmark      &
        \cmark                                                   & \xmark                          & \cmark                                              \\
        Prodigy~\cite{Mishchenko2024Prodigy}                     & \xmark                          & \xmark      & \xmark           & \cmark    &
        \xmark                                                                                                                                           \\
        D-Adaptation~\cite{Defazio2023DAdaptation}               & \xmark                          & \xmark      & \xmark           & \cmark
                                                                 & \xmark                                                                                \\
        UGM~\cite{Nesterov2015UniversalGM}                       & \cmark\tnote{{\color{blue}(*)}} & \cmark      &
        \cmark                                                   & \xmark                          & \xmark                                              \\
        DADA (\textbf{Ours})                                     & \cmark                          & \cmark      & \cmark           & \cmark    & \xmark \\
        \bottomrule
      \end{tabular}
      \begin{tablenotes}
        \footnotesize
        \item[{\cmark\color{blue}(*)}] Note that UGM uses a different definition
        of universality. They call their method universal because it works for
        H\"older-smooth functions, which are only a subset of the functions we
        consider.
      \end{tablenotes}
      \vspace{0.5em}
      \caption{A comparison of different adaptive algorithms to solve
        \eqref{eq:problem}. ``Universal'' means the algorithm achieves
        problem-dependent convergence rates across convex function
        classes by exploiting the local growth of the
        objective near the minimizer. ``Constraints'' means the
        algorithm can be applied to constrained problems.
        ``Unbounded domain'' means the algorithm can be applied to
        problems with unbounded feasible sets.
        ``Stochastic'' indicates that the algorithm is analyzed
        in the stochastic setting. ``No search'' means the algorithm does not use
        an internal search procedure.}
      \label{table:comparison-table}
    \end{threeparttable}
  }
\end{table}

The development of parameter-free first-order methods has received increasing
attention in both optimization and machine learning. A central goal in this
line of work is to design algorithms whose performance does not depend on prior
knowledge of problem's specific parameters, such as smoothness constants,
Lipschitz parameters, or distance to the minimizer—quantities that are rarely
known in practice.

Classical approaches to removing stepsize tuning include techniques such as
Polyak's stepsize rule~\cite{polyak1987} and doubling
schedules~\cite{streeter2012}. While effective in certain settings, these
strategies either rely on access to the optimal value or introduce additional
overhead through repeated restarts. In contrast, more recent parameter-free
methods aim to achieve near-optimal performance without requiring such
auxiliary procedures.

A large group of recent parameter-free methods is based on AdaGrad-type
conditioning~\cite{Duchi2011AdaGrad}. These methods adaptively accumulate
squared gradient norms to adjust the effective stepsize. This idea underlies
several recent distance-adaptation algorithms, including
DoG~\cite{Ivgi2023DoGIS}, DoWG~\cite{Khaled2023DoWG},
D-Adaptation~\cite{Defazio2023DAdaptation}, and
Prodigy~\cite{Mishchenko2024Prodigy}. Although these algorithms achieve
parameter-free convergence guarantees for nonsmooth Lipschitz or
Lipschitz-smooth objectives, their theoretical rates do not automatically
adapt to broader families of convex functions.
We summarize the main properties of the algorithms
we compare against in~\cref{table:comparison-table}.

Beyond AdaGrad-type schemes, coin-betting algorithms~\cite{orabona2016} provide
adaptive guarantees in online and stochastic optimization by treating learning
as a sequential investment game. In a different direction, Carmon and
Hinder~\cite{Carmon2022MakingSGDParameterFree} propose a bisection-based SGD
routine that adapts to the unknown smoothness or distance-to-optimum by
iteratively solving simpler subproblems. Both coin-betting and bisection
approaches are orthogonal to ours but share the goal of eliminating learning
rate tuning through adaptation mechanisms.

Another universal method worth noting is Nesterov's Universal Gradient Method
(UGM)~\cite{Nesterov2015UniversalGM}, which achieves optimal rates for
H\"older-smooth functions via adaptive line search. While UGM is often
described as “universal,” its scope is limited to smoothness-varying settings
and does not extend to broader function classes such as quasi-self-concordant
or high-order smooth functions. Moreover, its reliance on internal line search
procedures makes it less practical in constrained or composite problems.

\paragraph{Notation.}

In this text, we work in the space $\R^d$ equipped with the standard inner
product $\inp{\cdot}{\cdot}$ and the general Euclidean norm
$\norm{x} \defeq \inp{B x}{x}^{1 / 2}$,
where $B$ is a fixed symmetric positive definite matrix. The corresponding dual
norm is
defined in the standard way as
$\dnorm{s} \defeq \max_{\norm{x} = 1} \inp{s}{x} = \inp{s}{B^{-1}s}^{1 / 2}$.
Thus, for any $s, x \in \R^d$, we have the Cauchy-Schwarz inequality
$\abs{\inp{s}{x}} \le \dnorm{s}\norm{x}$.
The Euclidean ball of radius $r > 0$ centered at $x \in \R^d$ is defined as
$\Ball{x}{r} \defeq \SetBuilder{y \in \R^d}{\norm{y - x} \leq r}$.
For a convex function $\Map{f}{\R^d}{\RealFieldPlusInfty}$, we denote its
effective domain as
$\EffectiveDomain f \defeq \SetBuilder{x \in \R^d}{f(x) < + \infty}$.
The subdifferential of $f$ at a point
$x \in \EffectiveDomain f$ is denoted by $\Subdifferential f(x)$,
and $\nabla f(x) \in \Subdifferential f(x)$ denotes a subgradient.
We use $\nabla f(x)$ rather than $g \in \Subdifferential f(x)$
throughout the paper to keep the notation lightweight.

  \section{DADA Method}\label{sec:preliminaries}

\paragraph{Measuring the quality of solution.}

Given an approximate solution $x \in Q$ to problem~\eqref{eq:problem}
and an arbitrary subgradient $\nabla f(x) \in \Subdifferential f(x)$,
we measure the suboptimality of $x$ by the distance from $x^*$ to the
hyperplane
$\SetBuilder{y}{\inp{\nabla f(x)}{x - y} = 0}$:
\begin{align}\label{eq:v-definition}
  v(x) \DefinedEqual \frac{\inp{\nabla f(x)}{x - x^*}}{\dnorm{\nabla f(x)}}
  \quad
  \left(\geq 0\right).
\end{align}
This objective is meaningful because minimizing $v(x)$ also reduces the
corresponding function residual $f(x) - f^*$.
Indeed, there exists the following simple relationship between
$v(x)$ and the function residual
\cite[Section~3.2.2]{Nesterov2018LecturesOnCVX}
(see also \cref{lem:omega-upper} for the short proof):
\begin{align}\label{eq:omega-inequality}
  f(x) - f^* \le \omega(v(x)),
\end{align}
where
\begin{align}\label{eq:omega-definition}
  \omega(t)
  \DefinedEqual
  \max_{x \in \Ball{x^*}{t}} f(x) - f^*
\end{align}
measures the local growth of $f$ around the
solution $x^*$.
Note that inequality~\eqref{eq:omega-inequality} is nontrivial only when
$\Ball{x^*}{v(x)} \subseteq \EffectiveDomain f$.

By bounding $\omega(t)$, we can derive convergence-rate estimates that
simultaneously apply to a broad range of problem classes (we discuss several
  examples in \cref{sec:applications}).

\paragraph{The method.}

\begin{algorithm}[tb]
  \caption{General Scheme of DA}
  \label{alg:da}
  \begin{algorithmic}
    \Require
    $x_0 \in Q$, number of iterations $T \geq 1$,
    coefficients $(a_k)_{k = 0}^{T - 1}$,
    $(\beta_k)_{k = 1}^T$ with nondecreasing $\beta_k$
    \For{$k = 0, \ldots, T - 1$}
      \State Compute arbitrary $g_{k} \in \Subdifferential f(x_{k})$
      \State $x_{k + 1} = \argmin_{x \in Q} \bigl\{\sum_{i = 0}^{k}
          a_i \inp{g_i}{x - x_i} + \frac{\beta_{k + 1}}{2} \norm{x - x_0}^2 \bigr\}$
    \EndFor
    \Ensure $x^*_{T} = \argmin_{x \in \{x_0, \ldots, x_{T - 1}\}} f(x)$
  \end{algorithmic}
\end{algorithm}

Our algorithm is based on the general scheme of
DA~\cite{Nesterov2005PrimaldualSM} shown in \cref{alg:da}.
Using a standard (sub)gradient method with time-varying coefficients
is also possible but requires either short steps by fixing the number of
iterations in advance, or paying an extra logarithmic factor in the convergence
rate~\cite[Section~3.2.3]{Nesterov2018LecturesOnCVX}.

The classical method of Weighted DA (WDA) selects the coefficients
$a_k = \frac{\hat{D}_0}{\dnorm{g_k}}$ and
$\beta_k = \BigTheta(\sqrt{k})$,
where $\hat{D}_0$ is a user-defined estimate of $D_0$.
The convergence is guaranteed for any value of $\hat{D}_0$ but one must
pay a multiplicative cost of $\rho^2$, where
$\rho \DefinedEqual \max\{\frac{\hat{D}_0}{D_0}, \frac{D_0}{\hat{D}_0}\}$,
if the parameter $D_0$ is unknown. This cost can be
significantly high if $D_0$ is not known almost exactly.
To address this issue, we propose DADA, which reduces the cost to a logarithmic
term, $\log^2 \rho$, offering a substantial improvement.

Specifically, our approach utilizes the following coefficients:
\begin{align}\label{eq:dada-stepsizes}
  \boxed{\,a_k = \frac{\Bar{r}_k}{\dnorm{g_k}}, \quad \beta_k = c \sqrt{k +
      1}\,},
  \quad
  \Bar{r}_k \DefinedEqual \max \{\max_{1 \leq t \le k} r_t, \Bar{r} \},
  \quad
  r_t \DefinedEqual \norm{x_0 - x_t},
\end{align}
where $\Bar{r} > 0$ is a parameter
and $c$ is a certain constant to be specified later. In what follows,
we assume w.l.o.g.\ that
$g_k \neq 0$ for all $0 \leq k \leq T - 1$ since otherwise the exact
solution has been found, and the method could be successfully terminated
before making $T$ iterations.

As we can see, the main difference between WDA and DADA, is that the latter
dynamically adjusts its estimate of~$D_0$ by exploiting $r_t$, the distance
between $x_t$ and the initial point $x_0$. This idea has been
explored in recent works~\cite{Carmon2022MakingSGDParameterFree,Ivgi2023DoGIS},
which similarly utilize $r_t$ in various ways. Other methods also attempt to
estimate this quantity using alternative strategies, based on DA
and the similar principle of employing an increasing sequence of lower bounds
for $D_0$
\cite{Defazio2023DAdaptation, Mishchenko2024Prodigy}.

The convergence guarantees for our method are provided in the result below:
\begin{theorem}\label{thm:main-theorem}
  Consider \cref{alg:da} for solving problem~\eqref{eq:problem} using the
  coefficients from \cref{eq:dada-stepsizes} with $c > \sqrt{2}$. Then,
  for any $T \geq 1$ and
  $v_T^* \DefinedEqual \min_{0 \leq k \leq T - 1} v(x_k)$,
  it holds that
  \begin{align*}
    f(x_T^*) - f^* \le \omega(v_T^*),
  \end{align*}
  and
  \begin{equation}
    \label{eq:main-theorem:convergence-rate}
    v_T^* \leq \frac{e D}{\sqrt{T}} \log \frac{e \bar{D}}{\bar{r}},
  \end{equation}
  where $\bar{D} \DefinedEqual \max\{\bar{r}, \frac{2c}{c - \sqrt{2}} D_0\}$
  and $D \DefinedEqual \sqrt{2} (c D_0 + \frac{1}{c} \bar{D})$.
  Consequently, for a given $\delta > 0$, it holds that $v_T^* \leq \delta$
  whenever $T \geq T_v(\delta)$, where
  \begin{align*}
    T_v(\delta) \DefinedEqual
    \frac{e^2 D^2}{\delta^2} \log^2 \frac{e \bar{D}}{\bar{r}}.
  \end{align*}
\end{theorem}
Let us provide a proof sketch for \cref{thm:main-theorem} here and
defer the detailed proof to \cref{sec:main-theorem-v}.
We begin by applying the standard result for DA (\cref{lem:da-convergence}),
which holds for any choice of coefficients $a_k$ and $\beta_k$:
\begin{align*}
  \sum_{i = 0}^{k - 1} a_i v_i \dnorm{g_i} +
  \frac{\beta_{k}}{2}D_k^2 \le \frac{\beta_k}{2}D_0^2 + \sum_{i =
    0}^{k - 1}\frac{a_i^2}{2\beta_i}\dnorm{g_i}^2,
\end{align*}
where $D_i = \norm{x_i - x^*}$ and $v_i = v(x_i)$ for all $i \geq 0$.
Use the specific choices for $a_k$ and $\beta_k$ as
defined in \cref{eq:dada-stepsizes}, we obtain
(see \cref{lem:coefficient-placing}):
\begin{equation}
  \label{eq:dada-convergence}
  \sum_{i = 0}^{k - 1} \Bar{r}_iv_i + \frac{c \sqrt{k + 1}}{2} D_k^2
  \leq
  \frac{c\sqrt{k + 1}}{2} D_0^2 + \frac{\sqrt{k}}{c} \Bar{r}_{k - 1}^2.
\end{equation}
Dropping the nonnegative $\Bar{r}_iv_i$ from the left-hand side, we can
show
by induction that $\Bar{r}_k$ is uniformly bounded (see \cref{lem:r-upper-d}):
\begin{align*}
  \Bar{r}_k \leq \bar{D},
\end{align*}
where $\bar{D}$ is the constant from \cref{thm:main-theorem}.
This bound is crucial to our analysis, as we need to eliminate $\Bar{r}_{k -
    1}$ from
the right-hand side of \cref{eq:dada-convergence}. Achieving this requires
selecting the coefficients precisely as defined in \cref{eq:dada-stepsizes},
which is
the primary difference compared to the standard DA method
\cite{Nesterov2005PrimaldualSM}. Next, using the inequality
$D_0^2 - D_k^2 \leq 2 r_k D_0$,
we get
\begin{align*}
  \sum_{i = 0}^{k - 1} \Bar{r}_i v_i
  \leq
  c \sqrt{k + 1} r_k D_0
  +
  \frac{\sqrt{k}}{c} \bar{r}_{k - 1}^2
  \leq
  \Bigl( c D_0 + \frac{1}{c} \bar{D} \Bigr) \Bar{r}_{k}
  \sqrt{k + 1}.
\end{align*}
After establishing this, the rest of the proof follows straightforwardly by
dividing both sides by $\sum_{i = 0}^{k - 1} \Bar{r}_i$ and applying the
following inequality (valid for any nondecreasing sequence $\bar{r}_k$, see
  \cref{lem:seq-min-sum}):
\begin{align*}
  \min_{1 \le k \le T} \frac{\Bar{r}_k}{\sum_{i = 0}^{k - 1} \Bar{r}_i} \leq
  \frac{(\frac{\Bar{r}_T}{\Bar{r}_0})^{\frac{1}{T}} \log\frac{e
    \Bar{r}_T}{\Bar{r}_0}}{T}.
\end{align*}
This gives us
\[
  v_T^*
  \leq
  \frac{D}{\sqrt{T}}
  \biggl( \frac{\bar{D}}{\bar{r}} \biggr)^{\frac{1}{T}}
  \log \frac{e \bar{D}}{\bar{r}},
\]
which is almost \cref{eq:main-theorem:convergence-rate} except for the extra
factor of $(\frac{\bar{D}}{\bar{r}})^{\frac{1}{T}}$.
This extra factor, however, is rather weak as it can be upper bounded by a
constant (say, $e \equiv \exp(1)$) whenever
$T \geq \log \frac{\bar{D}}{\bar{r}}$.
The case of $T \leq \log \frac{\bar{D}}{\bar{r}}$ is not interesting since
then \cref{eq:main-theorem:convergence-rate} holds trivially because,
for any $k \geq 0$, in view of \cref{eq:v-definition,lem:r-upper-d},
we have $v_k \leq D_k \leq D$.
According to \cref{thm:main-theorem}, our method converges for any
$c > \sqrt{2}$. To obtain the smallest complexity bound
(up to logarithmic factors), the value that minimizes this bound is
$c = 2 \sqrt{2}$.
A more detailed discussion of this choice is provided in
\cref{sec:discussion-on-c}.

  \section{Universality of DADA: Examples of
  Applications}\label{sec:applications}

Let us demonstrate that our method is \emph{universal} in the sense that it
simultaneously works for multiple problem classes without the need
for choosing different parameters for each of these function classes. For
simplicity, we assume that $\nabla f(x^*) = 0$ (this happens, in particular,
  when our problem~\eqref{eq:problem} is unconstrained) and measure the
$\epsilon$-accuracy in terms of the function residual.
This assumption is made only to keep the discussion clean and readable,
and it also reflects an important practical setting (unconstrained problems,
  or constrained problems with $x^*$ in the interior of $Q$).
The general constrained case, where $\nabla f(x^*) \neq 0$, is covered in
\cref{sec:convergence-proof-examples}; there, for some function classes
(e.g., Lipschitz-smooth), the complexity bounds include an additional
$\dnorm{\nabla f(x^*)}^2 / \epsilon^2$ term that may affect the rate.
We also assume, for simplicity, that the objective function satisfies
all necessary inequalities on the entire space, but all our results
still hold if they are satisfied only
locally at $x^*$ (see \cref{sec:convergence-proof-examples}).
To simplify the notation, we also denote
$\log_+ t \DefinedEqual 1 + \log t$ and
$\bar{D}_{0} \DefinedEqual \max\{ \bar{r}, \norm{x_0 - x^*} \}$,
where $\bar{r}$ is the parameter of our method.

\subsection{Nonsmooth Lipschitz Functions}\label{sec:lipschitz-function}

\begin{assumption}\label{assumption:locally-lipschitz-function}
  The function $f$ in problem~\eqref{eq:problem} is locally Lipschitz at $x^*$.
  Specifically, for any $x \in \Ball{x^*}{\rho}$,
  the following inequality holds:
  \begin{align}\label{eq:locally-lipschitz-function}
    f(x) - f^* \leq L_0 \norm{x - x^*},
  \end{align}
  where $L_0, \rho > 0$ are fixed constants.
\end{assumption}

\begin{lemma}\label{lem:lipschitz-omega}
  Let $f$ be locally $L_0$-Lipschitz at $x^*$
  (\cref{assumption:locally-lipschitz-function}).
  Then, $\omega(t) \leq \epsilon$ for any given $\epsilon > 0$ whenever
  $t \le \delta(\epsilon)$, where
  \begin{align*}
    \delta(\epsilon) \DefinedEqual \min\left\{\frac{\epsilon}{L_0}, \rho\right\}.
  \end{align*}
\end{lemma}

\begin{proof}
  According to \eqref{eq:locally-lipschitz-function},
  for any $0 \le t \le \rho$, we have
  \begin{align*}
    \omega(t) \leq L_0 t.
  \end{align*}
  Making the right-hand side $\leq \epsilon$, we get the claim.
\end{proof}

Combining \cref{thm:main-theorem,lem:lipschitz-omega}, we get the following
complexity result.

\begin{corollary}\label{thm:nonsmooth-lipschitz-convergence}
  Consider problem~\eqref{eq:problem}
  under~\cref{assumption:locally-lipschitz-function}.
  Let \cref{alg:da} with coefficients \eqref{eq:dada-stepsizes}
  be applied for solving this problem.
  Then, $f(x_T^*) - f^* \leq \epsilon$ for any given $\epsilon > 0$ whenever
  $T \geq T(\epsilon)$, where
  \begin{align*}
    T(\epsilon)
    = \max\left\{\frac{L_0^2}{\epsilon^2}, \frac{1}{\rho^2}\right\}
    e^2 D^2 \log^2 \frac{e \bar{D}}{\Bar{r}},
  \end{align*}
  and the constants $D$ and $\bar{D}$
  are as defined in \cref{thm:main-theorem}.
\end{corollary}

\subsection{Lipschitz-Smooth Functions}\label{sec:lip-grad-func}

\begin{assumption}\label{assumption:locally-lipschitz-smooth-function}
  The function $f$ in problem~\eqref{eq:problem} is
  locally Lipschitz-smooth at $x^*$. Specifically, for any
  $x \in \Ball{x^*}{\rho}$, the following inequality holds:
  \begin{align}\label{eq:locally-lipschitz-smooth-function}
    f(x) \leq f^* + \inp{\nabla f(x^*)}{x - x^*} +
    \frac{L_1}{2} \norm{x - x^*}^2,
  \end{align}
  where $L_1, \rho > 0$ are fixed constants.
\end{assumption}

\begin{lemma}\label{lem:smooth-omega}
  Assume that $f$ is locally Lipschitz-smooth at $x^*$ with constant $L_1$
  (\cref{assumption:locally-lipschitz-smooth-function}).
  Then, $\omega(t) \leq \epsilon$ for any given $\epsilon > 0$ whenever
  $t \le \delta(\epsilon)$, where
  \begin{align*}
    \delta(\epsilon) \DefinedEqual \min\left\{
      \sqrt{\frac{\epsilon}{L_1}},
      \frac{\epsilon}{2\dnorm{\nabla f(x^*)}},
      \rho
    \right\}.
  \end{align*}
\end{lemma}
\begin{proof}
  According to \eqref{eq:locally-lipschitz-smooth-function},
  for any $x \in \Ball{x^*}{\rho}$, we have
  \begin{align*}
    f(x) - f^*
    \leq \dnorm{\nabla f(x^*)}\norm{x - x^*} + \frac{L_1}{2} \norm{x - x^*}^2.
  \end{align*}
  Hence, for any $0 \le t \le \rho$,
  \begin{align*}
    \omega(t) & \le \frac{L_1}{2} t^2 + \dnorm{\nabla f(x^*)} t.
  \end{align*}
  To make the right-hand side $\leq \epsilon$, it suffices to ensure that each of
  the two terms is $\leq \frac{\epsilon}{2}$:
  \begin{align*}
    \frac{L_1}{2} t^2 \leq \frac{\epsilon}{2},
    \qquad
    \dnorm{\nabla f(x^*)} t \leq \frac{\epsilon}{2}.
  \end{align*}
  Solving this system of inequalities, we get the claim.
\end{proof}

Combining \cref{thm:main-theorem,lem:smooth-omega}, we get the following
complexity result.

\begin{corollary}\label{thm:lipschitz-smooth-convergence}
  Consider problem~\eqref{eq:problem} under
  \cref{assumption:locally-lipschitz-smooth-function}.
  Let \cref{alg:da} with coefficients \eqref{eq:dada-stepsizes}
  be applied for solving this problem.
  Then, $f(x_T^*) - f^* \leq \epsilon$ for any given $\epsilon > 0$ whenever
  $T \geq T(\epsilon)$, where
  \begin{align*}
    T(\epsilon)
    = \max \left\{\frac{L_1}{\epsilon},
      \frac{4 \dnorm{\nabla f(x^*)}^2}{\epsilon^2}, \frac{1}{\rho^2}\right\}
    e^2 D^2 \log^2 \frac{e \bar{D}}{\Bar{r}},
  \end{align*}
  and the constants $D$ and $\bar{D}$
  are as defined in \cref{thm:main-theorem}.
\end{corollary}

\subsection{H\"older-Smooth Functions}\label{sec:holder-smooth}

\begin{assumption}\label{assumption:locally-holder-smooth-function}
  The function $f$ in problem~\eqref{eq:problem} is
  locally H\"older-smooth at $x^*$.
  Specifically, for any $x \in \Ball{x^*}{\rho}$,
  the following inequality holds:
  \begin{align}\label{eq:locally-holder-smooth-function}
    f(x) \leq f^* + \inp{\nabla f(x^*)}{x - x^*} +
    \frac{H_{\nu}}{1 + \nu} \norm{x - x^*}^{1 + \nu},
  \end{align}
  where $\nu \in [0, 1]$ and $H_{\nu}, \rho > 0$ are fixed constants.
\end{assumption}

\begin{lemma}\label{lem:holder-smooth-omega}
  Let $f$ be locally $(\nu, H_{\nu})$-H\"older-smooth at $x^*$
  (\cref{assumption:locally-holder-smooth-function}).
  Then, $\omega(t) \leq \epsilon$ for any given $\epsilon > 0$ whenever
  $t \le \delta(\epsilon)$, where
  \begin{align*}
    \delta(\epsilon) \DefinedEqual \min \left\{
      \left[\frac{(1 + \nu) \epsilon}{2 H_{\nu}}\right]^{\frac{1}{1 + \nu}},
      \frac{\epsilon}{2\dnorm{\nabla f(x^*)}},
      \rho
    \right\}.
  \end{align*}
\end{lemma}
\begin{proof}
  According to \eqref{eq:locally-holder-smooth-function},
  for any $x \in \Ball{x^*}{\rho}$, we have
  \begin{align*}
    f(x) - f^*
    \leq \dnorm{\nabla f(x^*)} \norm{x - x^*} +
    \frac{H_{\nu}}{1 + \nu} \norm{x - x^*}^{1 + \nu}.
  \end{align*}
  Hence, for any $0 \le t \le \rho$,
  \begin{align*}
    \omega(t) \leq \dnorm{\nabla f(x^*)} t
    + \frac{H_{\nu}}{1 + \nu} t^{1 + \nu}.
  \end{align*}
  To make the right-hand side of the last inequality $\leq \epsilon$,
  it suffices to ensure that each of the two terms is $\leq \frac{\epsilon}{2}$:
  \begin{align*}
    \dnorm{\nabla f(x^*)} t \leq \frac{\epsilon}{2},
    \qquad
    \frac{H_{\nu}}{1 + \nu} t^{1 + \nu} \leq \frac{\epsilon}{2}.
  \end{align*}
  Solving this system of inequalities, we get the claim.
\end{proof}

Combining \cref{thm:main-theorem,lem:holder-smooth-omega}, we get the following
complexity result.

\begin{corollary}\label{thm:holder-smooth-convergence}
  Consider problem~\eqref{eq:problem} under
  \cref{assumption:locally-holder-smooth-function}.
  Let \cref{alg:da} with coefficients \eqref{eq:dada-stepsizes}
  be applied for solving this problem.
  Then, $f(x_T^*) - f^* \leq \epsilon$ for any given $\epsilon > 0$ whenever
  $T \geq T(\epsilon)$, where
  \begin{align*}
    T(\epsilon) = \max \left\{
      \left[\frac{2 H_{\nu}}{(1 + \nu) \epsilon}\right]^{\frac{2}{1 + \nu}},
      \frac{4 \dnorm{\nabla f(x^*)}^2}{\epsilon^2}, \frac{1}{\rho^2} \right\}
    e^2 D^2 \log^2 \frac{e \bar{D}}{\Bar{r}},
  \end{align*}
  and the constants $D$ and $\bar{D}$
  are as defined in \cref{thm:main-theorem}.
\end{corollary}

\subsection{Functions with Lipschitz High-Order Derivative}
\label{sec:pth-lip-func}

\begin{assumption}\label{assumption:locally-pth-grad-lipschitz-function}
  The function $f$ in problem~\eqref{eq:problem} is such that
  its $p$th derivative is locally $L_p$-Lipschitz at $x^*$.
  Specifically, $f$ is $p$ times differentiable on $\Ball{x^*}{\rho}$,
  and, for any $x \in \Ball{x^*}{\rho}$, the following inequality holds:
  \begin{align*}
    \norm{\nabla^p f(x) - \nabla^p f(x^*)} \le L_p \norm{x - x^*}.
  \end{align*}
  where $L_p, \rho > 0$ are fixed constants.
\end{assumption}

The \cref{assumption:locally-pth-grad-lipschitz-function}
immediately implies the following global upper bound on the function value:
\begin{align}\label{eq:locally-pth-grad-lipschitz-function}
  f(x) \leq
  f^* + \sum_{i = 1}^{p}
  \frac{1}{i!} \nabla^i f(x^*)[x - x^*]^i
  + \frac{L_p}{(p + 1)!} \norm{x - x^*}^{p + 1}.
\end{align}

\begin{lemma}\label{lem:pth-grad-lipschitz-omega}
  Assume that $f$ has locally $L_p$-Lipschitz $p$th derivative at $x^*$
  (\cref{assumption:locally-pth-grad-lipschitz-function}).
  Then, $\omega(t) \leq \epsilon$ for any given $\epsilon > 0$ whenever
  $t \leq \delta(\epsilon)$, where
  \begin{align*}
    \delta(\epsilon) \DefinedEqual \min \left\{
      \min_{2 \leq i \leq p}
      \left[
        \frac{i! \, \epsilon}{(p + 1) \norm{\nabla^i f(x^*)}}\right]^{\frac{1}{i}},
      \left[\frac{p! \,\epsilon}{L_p}\right]^{\frac{1}{p + 1}},
      \frac{\epsilon}{(p + 1) \dnorm{\nabla f(x^*)}},
      \rho
    \right\}.
  \end{align*}
\end{lemma}
\begin{proof}
  According to \eqref{eq:locally-pth-grad-lipschitz-function},
  for any $x \in \Ball{x^*}{\rho}$, we have
  \begin{align*}
    f(x) - f^*
    \leq \dnorm{\nabla f(x^*)} \norm{x - x^*} +
    \sum_{i = 2}^{p} \frac{1}{i!} \norm{\nabla^i f(x^*)} \norm{x - x^*}^i +
    \frac{L_p}{(p + 1)!} \norm{x - x^*}^{p + 1}.
  \end{align*}
  Therefore, for any $0 \le t \le \rho$, we have
  \begin{align*}
    \omega(t) & \leq
    \dnorm{\nabla f(x^*)} t +
    \sum_{i = 2}^{p} \frac{1}{i!} \norm{\nabla^i f(x^*)} t^i +
    \frac{L_p}{(p + 1)!} t^{p + 1}.
  \end{align*}
  To make the right-hand side $\leq \epsilon$, it suffices to ensure that each of
  the following inequalities holds:
  \begin{align*}
    \dnorm{\nabla f(x^*)} t \leq \frac{\epsilon}{p + 1}, \quad
    \frac{1}{i!} \norm{\nabla^i f(x^*)} t^i \leq
    \frac{\epsilon}{p + 1}, \quad \frac{L_p}{(p + 1)!} t^{p + 1} \leq
    \frac{\epsilon}{p + 1}, \quad i = 2, \ldots, p.
  \end{align*}
  Solving this system of inequalities, we get the claim.
\end{proof}

Combining \cref{thm:main-theorem,lem:pth-grad-lipschitz-omega}, we get the
following complexity result.

\begin{corollary}\label{thm:high-order-derivative-lipschitz-convergence}
  Consider problem~\eqref{eq:problem} under
  \cref{assumption:locally-pth-grad-lipschitz-function}.
  Let \cref{alg:da} with coefficients \eqref{eq:dada-stepsizes}
  be applied for solving this problem.
  Then, $f(x_T^*) - f^* \leq \epsilon$ for any given $\epsilon > 0$ whenever
  $T \geq T(\epsilon)$, where
  \begin{align*}
    T(\epsilon) =
    \max \Biggl\{
      & \max_{2 \leq i \leq p} \left[\frac{(p + 1)
          \norm{\nabla^i f(x^*)}}{i! \, \epsilon}\right]^{\frac{2}{i}}, \\
      & \left[\frac{L_p}{p! \,\epsilon}\right]^{\frac{2}{p + 1}},
      \frac{(p + 1)^2 \dnorm{\nabla f(x^*)}^2}{\epsilon^2},
      \frac{1}{\rho^2}
    \Biggr\} e^2 D^2 \log^2 \frac{e \bar{D}}{\Bar{r}},
  \end{align*}
  and the constants $D$ and $\bar{D}$
  are as defined in \cref{thm:main-theorem}.
\end{corollary}

\subsection{Quasi-Self-Concordant Functions}

\begin{assumption}\label{assumption:locally-qsc-function}
  The function $f$ in problem~\eqref{eq:problem} is
  Quasi-Self-Concordant (QSC) in a neighborhood of $x^*$.
  Specifically, it is three times differentiable in a neighborhood of $x^*$
  and for any $x \in \Ball{x^*}{\rho}$ and arbitrary directions
  $u, v \in \R^d$, the following inequality holds:
  \begin{align*}
    \nabla^3 f(x)[u, u, v] \le M \inp{\nabla^2 f(x)u}{u} \norm{v},
  \end{align*}
  where $M \geq 0$ and $\rho > 0$ are fixed constants.
\end{assumption}

The following lemma provides an important global upper bound on the
function value for QSC functions.
\begin{lemma}{\cite[Lemma~2.7]{Doikov2023QuasiSelfConcordant}}
  \label{lem:quasi-self-concordant-descent}
  Let f be QSC with the parameter $M$.
  Then, for any $x, y \in \EffectiveDomain f$, the following inequality holds:
  \begin{align*}
    f(y) \le f(x) + \inp{\nabla f(x)}{y - x} + \inp{\nabla^2 f(x)(y - x)}{y - x}
    \varphi(M\norm{y - x}),
  \end{align*}
  where $\varphi(t) \defeq \frac{e^t - t - 1}{t^2}$.
\end{lemma}

\begin{lemma}\label{lem:qsc-omega}
  Assume that $f$ is a locally QSC function at $x^*$ with constant $M$
  (\cref{assumption:locally-qsc-function}).
  Then, $\omega(t) \leq \epsilon$ for any given $\epsilon > 0$ whenever
  $t \leq \delta(\epsilon)$, where
  \begin{align*}
    \delta(\epsilon) \DefinedEqual \min \left\{
      \frac{1}{M},
      \sqrt{\frac{\epsilon}{2 (e - 2) \norm{\nabla^2 f(x^*)}}},
      \frac{\epsilon}{2 \dnorm{\nabla f(x^*)}},
      \rho
    \right\}.
  \end{align*}
\end{lemma}
\begin{proof}
  According to \cref{lem:quasi-self-concordant-descent},
  for any $x \in \Ball{x^*}{\rho}$, we
  have
  \begin{align*}
    f(x) - f^* & \leq \inp{\nabla f(x^*)}{x - x^*} + \inp{\nabla^2 f(x^*)(x -
        x^*)}{x - x^*} \varphi(M\norm{x - x^*}) \\
    & \le \dnorm{\nabla f(x^*)} \norm{x - x^*} + \norm{\nabla^2 f(x^*)} \norm{x -
      x^*}^2 \varphi(M\norm{x - x^*}).
  \end{align*}
  Therefore, for any $0 \le t \le \rho$, we get
  \begin{align}\label{eq:omega-qsc}
    \omega(t)
    \leq \dnorm{\nabla f(x^*)} t + \norm{\nabla^2 f(x^*)} t^2 \varphi(M t),
  \end{align}
  where we have used the fact that $\varphi(\cdot)$ is an increasing function.

  Note that, for any $0 \leq t \leq \frac{1}{M}$, we can estimate
  $\varphi(M t) \leq \varphi(1) = e - 2$. Substituting this bound into
  \eqref{eq:omega-qsc}, we obtain
  \begin{align*}
    \omega(t) \leq
    \dnorm{\nabla f(x^*)} t + (e - 2) \norm{\nabla^2 f(x^*)} t^2.
  \end{align*}
  To make the right-hand side $\leq \epsilon$, it suffices to ensure that each of
  the two terms is $\leq \frac{\epsilon}{2}$:
  \begin{align*}
    \dnorm{\nabla f(x^*)} t \leq \frac{\epsilon}{2},
    \qquad
    (e - 2) \norm{\nabla^2 f(x^*)} t^2 \leq \frac{\epsilon}{2}.
  \end{align*}
  Solving this system of inequalities, we get the claim.
\end{proof}

Combining \cref{thm:main-theorem,lem:qsc-omega}, we get the following
complexity result.

\begin{corollary}\label{thm:qsc-convergence}
  Consider problem~\eqref{eq:problem} under
  \cref{assumption:locally-qsc-function}.
  Let \cref{alg:da} with coefficients \eqref{eq:dada-stepsizes}
  be applied for solving this problem.
  Then, $f(x_T^*) - f^* \leq \epsilon$ for any given $\epsilon > 0$ whenever
  $T \geq T(\epsilon)$, where
  \begin{align*}
    T(\epsilon) =
    \max \left\{
      M^2,
      \frac{2 (e - 2) \norm{\nabla^2 f(x^*)}}{\epsilon},
      \frac{4 \dnorm{\nabla f(x^*)}^2}{\epsilon^2},
      \frac{1}{\rho^2}
    \right\}
    e^2 D^2
    \log^2 \frac{e \bar{D}}{\Bar{r}},
  \end{align*}
  and the constants $D$ and $\bar{D}$
  are as defined in \cref{thm:main-theorem}.
\end{corollary}

\subsection{\texorpdfstring{$(L_0, L_1)$}{(L0, L1)}-Smooth Functions}
\label{sec:l0-l1-smooth}

Let us now consider the case when $Q = \R^d$ and $f$ is $(L_0, L_1)$-smooth
\cite{Zhang2020Why},
meaning that for any $x \in \R^d$,
\begin{align*}
  \norm{\nabla^2 f(x)} \leq L_0 + L_1 \dnorm{\nabla f(x)},
\end{align*}
where $L_0, L_1 \geq 0$ are fixed constants.

\begin{lemma}{\cite[Lemma~2.2]{Vankov2024OptimizingL0L1SmoothFunctions}}
  \label{lem:l0-l1-descent}
  Let $f$ be $(L_0, L_1)$-smooth. Then, for any $x, y \in \R^d$, it holds that
  \begin{align*}
    f(y) \le f(x) + \inp{\nabla f(x)}{y - x} + \frac{L_0 + L_1 \dnorm{\nabla
        f(x)}}{L_1^2} \xi(L_1 \norm{y - x}),
  \end{align*}
  where $\xi(t) \DefinedEqual e^t - t - 1$.
\end{lemma}

\begin{lemma}\label{lem:l0-l1-omega}
  Assume that $f$ is an $(L_0, L_1)$-smooth function.
  Then, $\omega(t) \leq \epsilon$ for any given $\epsilon > 0$ whenever
  $t \leq \delta(\epsilon)$, where
  \begin{align*}
    \delta(\epsilon) \DefinedEqual \min \left\{
      \frac{1}{L_1},
      \sqrt{\frac{2 \epsilon}{3 (L_0 + L_1 \dnorm{\nabla f(x^*)})}},
      \frac{\epsilon}{2 \dnorm{\nabla f(x^*)}}
    \right\}.
  \end{align*}
\end{lemma}
\begin{proof}
  According to \cref{lem:l0-l1-descent}, for any $x \in \R^d$, we have
  \begin{align*}
    f(x) - f^* & \leq \inp{\nabla f(x^*)}{x - x^*} + \frac{L_0 + L_1
      \dnorm{\nabla f(x^*)}}{L_1^2} \xi(L_1 \norm{x - x^*}) \\
    & \leq \dnorm{\nabla f(x^*)} \norm{x - x^*} + \frac{L_0 + L_1
      \dnorm{\nabla f(x^*)}}{L_1^2} \xi(L_1 \norm{x - x^*})
  \end{align*}
  Therefore, for any $t \geq 0$, we get
  \begin{align}\label{eq:omega-l0-l1}
    \omega(t)
    \leq \dnorm{\nabla f(x^*)} t +
    \frac{L_0 + L_1 \dnorm{\nabla f(x^*)}}{L_1^2} \xi(L_1 t),
  \end{align}
  where the second inequality uses the fact that $\xi(x)$
  is an increasing function.

  Note that, for any $0 \leq t \leq \frac{1}{L_1}$, we can estimate
  \begin{align*}
    \xi(L_1 t) \leq \frac{L_1^2 t^2}{2 (1 - \frac{L_1 t}{3})}
    \leq \frac{3}{4} L_1^2 t^2.
  \end{align*}
  Substituting this bound into
  \eqref{eq:omega-l0-l1}, we obtain:
  \begin{align*}
    \omega(t)
    \leq
    \dnorm{\nabla f(x^*)} t +
    \frac{3(L_0 + L_1 \dnorm{\nabla f(x^*)})}{4} t^2.
  \end{align*}
  To make the right-hand side of the last inequality $\leq \epsilon$, it suffices
  to ensure that each of the two terms is $\leq \frac{\epsilon}{2}$:
  \begin{align*}
    \dnorm{\nabla f(x^*)} t \leq \frac{\epsilon}{2},
    \qquad
    \frac{3(L_0 + L_1 \dnorm{\nabla f(x^*)})}{4} t^2
    \leq \frac{\epsilon}{2}.
  \end{align*}
  Solving this system of inequalities, we get the claim.
\end{proof}

Combining \cref{thm:main-theorem,lem:l0-l1-omega}, we get the following
complexity result.

\begin{corollary}\label{thm:l0-l1-smooth-convergence}
  Consider problem~\eqref{eq:problem} under the assumption that
  $f$ is an $(L_0, L_1)$-smooth function.
  Let \cref{alg:da} with coefficients \eqref{eq:dada-stepsizes}
  be applied for solving this problem.
  Then, $f(x_T^*) - f^* \leq \epsilon$ for any given $\epsilon > 0$ whenever
  $T \geq T(\epsilon)$, where
  \begin{align*}
    T(\epsilon) =
    \max \left\{
      L_1^2,
      \frac{3 (L_0 + L_1 \dnorm{\nabla f(x^*)})}{2 \epsilon},
      \frac{4 \dnorm{\nabla f(x^*)}^2}{\epsilon^2}
    \right\}
    e^2 D^2
    \log^2 \frac{e \bar{D}}{\Bar{r}},
  \end{align*}
  and the constants $D$ and $\bar{D}$
  are as defined in \cref{thm:main-theorem}.
\end{corollary}

  \section{Experiments}\label{sec:experiments}

To evaluate the efficiency of our proposed method, DADA,
we conduct a series of experiments on convex optimization problems.
Our goal is to demonstrate the effectiveness of DADA in achieving competitive
performance across various function classes
\emph{without any hyperparameter tuning}.

We compare DADA against state-of-the-art distance-adaptation algorithms,
namely, DoG~\cite{Ivgi2023DoGIS} and Prodigy~\cite{Mishchenko2024Prodigy},
using their official implementations without any modifications.
We also consider the Universal Gradient Method (UGM)
from~\cite{Nesterov2015UniversalGM} and the classical
Weighted Dual Averaging (WDA) method~\cite{Nesterov2005PrimaldualSM}.
For UGM, we choose the initial value of the line-search parameter~$L_0 = 1$
and set the target accuracy to $\epsilon = 10^{-6}$.
For WDA, we use the coefficients $a_k = \frac{D_0}{\dnorm{g_k}}$ and
$\beta_k = \sqrt{k}$, where $D_0 = \Norm{x_0 - x^*}$.

For each method, we plot the best function value among all the test points
generated by the algorithm against the number of first-order oracle calls.
We set the starting point to $x_0 = (1, \ldots, 1)$ and select the initial
guess for the distance to the solution as
$\bar{r} = \delta (1 + \norm{x_0})$. This choice ensures a fair
comparison between DADA and DoG~\cite{Ivgi2023DoGIS}, as DoG employs a similar
initialization strategy. In all experiments, we fix $\delta = 10^{-6}$.
Additionally, we conduct a separate experiment to evaluate the sensitivity of
DADA to the choice of $\delta$.

We have several experiments on different problem classes. However,
due to space constraints, we present only a single representative experiment in
this section. The remaining experiments can be found in
\cref{sec:additional-experiments}.

\paragraph{Worst-case function.}

\begin{figure}[tb]
  \centering
  \includegraphics[width=\textwidth]{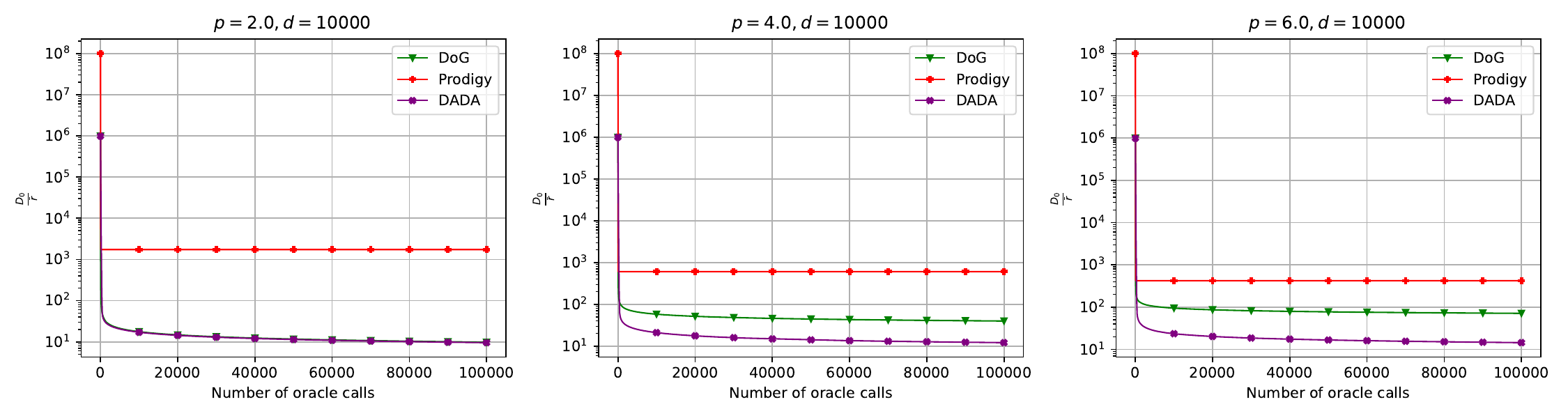}
  \caption{
    The ratio $\frac{D}{\bar{r}_t}$
    for the worst-case function with different optimal points $x^*$.
  }
  \label{fig:worst-instance-estimate-error}
\end{figure}

In addition to the experiments presented in \cref{sec:experiments},
we evaluate the estimation error of $D_0$
for Prodigy, DoG, and DADA throughout the optimization process,
as shown in \cref{fig:worst-instance-estimate-error}.
The figure illustrates that while Prodigy's estimate of $D_0$,
shows some improvement over time, it remains noticeably inaccurate.
Moreover, for DoG, the estimate deteriorates as $p$ increases, a behavior that
is not observed with DADA,
whose estimate remains stable across different values of $p$.

  \section{Discussion}\label{sec:conclusion}

\paragraph{Comparison with recent distance-adaptation methods.}

Let us briefly compare our method with several recently proposed parameter-free
algorithms, namely, DoG~\cite{Ivgi2023DoGIS}, DoWG~\cite{Khaled2023DoWG},
D-Adaptation~\cite{Defazio2023DAdaptation} and
Prodigy~\cite{Mishchenko2024Prodigy}.

To begin, we clarify the key differences between our method and other existing
gradient methods using the distance-adaptation technique.
One immediate difference is that we use DA instead
of the classical (sub)gradient method employed by DoG.
We could also instantiate our approach using the standard subgradient method
instead of DA.
However, doing so would either require fixing the number of iterations in
advance or would worsen the overall complexity
by an additional polylogarithmic factor in the target accuracy.
However, the most significant difference lies in how the sequence of gradients
is handled. In contrast to existing distance-adaptation methods,
which follow the AdaGrad~\cite{Duchi2011AdaGrad}
principle of accumulating squared gradient norms,
our method simply normalizes $g_k$ by its own norm.
This modification makes our method universal,
ensuring that $v_T^*$---the distance from $x^*$ to the supporting
hyperplane---converges to zero,
which is not known to be the case for DoG, even for deterministic problems.

Both DoG and DoWG employ a similar approach to estimate $D_0 = \norm{x_0 -
    x^*}$
and achieve comparable convergence rates for Lipschitz-smooth and nonsmooth
functions.
Similarly to our approach, DoWG considers only the
deterministic case, but with an additional assumption of the bounded
feasible set. They
have a different definition of universality, considering only Lipschitz-smooth
and nonsmooth settings.

Furthermore, to the best of our knowledge, these D-Adaptation and Prodigy
have not been extended to constrained optimization.
Nonetheless, their methods yield notable results in experiments,
demonstrating strong empirical performance.

It is important to emphasize that the advantage of our method over DoG does not
lie in guaranteeing convergence.
Indeed, \cite[Theorem 1]{Ivgi2023DoGIS} shows that DoG asymptotically converges
to a minimizer for any convex function, with a complexity of
$\BigOTilde\bigl( \frac{L_R^2 \bar{D}_0^2}{\epsilon^2} \bigr)$,
where $L_R = \max_{x \in \Ball{x_0}{R}} \dnorm{\nabla f(x)}$ for
$R = 3 \bar{D}_0$, and $\epsilon$ denotes the target accuracy in the function
value.
However, this complexity bound has a critical drawback---it remains inversely
proportional to $\epsilon^2$ across \emph{all} function classes, which is not
the case for our method.
For illustration, consider the setting where $f$ has an $L_2$-Lipschitz
Hessian.
Further, assume for simplicity that the problem is unconstrained and that
$\norm{\nabla^2 f(x^*)}$ is zero (or negligibly small).
In this case, the above complexity bound for DoG becomes\footnote{%
  Here, we use the fact that for functions with Lipschitz Hessian,
  $L_R = \BigO\left(L_1^* R + \frac{L_2}{2} R^2\right)$, where $L_1^* =
    \norm{\nabla^2 f(x^*)}$.
}
$\BigOTilde\bigl( \frac{L_2^2 \bar{D}_0^6}{\epsilon^2} \bigr)$,
which is substantially worse than
$
  \BigOTilde\bigl(
    \frac{L_2^{2 / 3} \bar{D}_0^2}{\epsilon^{2 / 3}}
  \bigr)
$
for DADA (see \cref{thm:high-order-derivative-lipschitz-convergence}).
Thus, in comparison to DoG, our method provides significantly stronger
efficiency guarantees and exhibits automatic acceleration under favorable
conditions for a considerably broader family of function classes.

\paragraph{Conclusion.}

We proposed DADA, a new adaptive and universal optimization method that
extends the classical Dual Averaging algorithm with a novel distance adaptation
mechanism.
Our method achieves competitive rates across a wide class of convex
problems---including Lipschitz, Lipschitz-smooth, H\"older-smooth,
quasi-self-concordant (QSC), and $(L_0, L_1)$-smooth functions---without
requiring parameter tuning or knowledge of smoothness constants.
In contrast to recent approaches such as DoG, DoWG, D-Adaptation, and Prodigy,
DADA seamlessly accommodates both constrained and unconstrained settings, and
does so without requiring restarts or line searches.

DADA provides a unified and adaptive framework for convex optimization with
convergence guarantees under minimal assumptions.
Future work includes extending DADA to stochastic and nonconvex optimization,
and evaluating its empirical performance in large-scale learning tasks.

  \printbibliography

  \newpage
  \appendix
  \section{Auxiliary Results}\label{sec:auxiliary-results}
The following result has been established in prior works such as
\cite[Lemma~30]{Liu2023StochasticNC}.
We include the proof here for the reader's convenience.
\begin{lemma}\label{lem:seq-min-sum}
  Let $(d_i)_{i = 0}^\infty$ be a positive nondecreasing sequence.
  Then for any $T \geq 1$,
  \begin{align*}
    \min_{1 \le k \le T} \frac{d_k}{\sum_{i = 0}^{k - 1} d_i} \le
    \frac{(\frac{d_{T}}{d_0})^{\frac{1}{T}} \log\frac{e
      d_{T}}{d_0}}{T}.
  \end{align*}
\end{lemma}

\begin{proof}
  Let $A_k \defeq \frac{1}{d_k} \sum_{i = 0}^{k - 1} d_i$ for each
  $k \geq 0$ (so that $A_0 = 0$). Then, for each $k \geq 0$, we have
  \begin{align*}
    d_{k + 1} A_{k + 1} - d_{k} A_{k} = d_{k},
  \end{align*}
  which implies that
  \begin{align*}
    \frac{d_{k}}{d_{k + 1}}
    = A_{k + 1} - \frac{d_{k}}{d_{k + 1}} A_{k}
    = A_{k + 1} - A_{k} + \left(1 - \frac{d_{k}}{d_{k + 1}}\right)A_{k}.
  \end{align*}
  Summing up these identities for all $0 \le k \le T - 1$, we get
  \begin{align*}
    S_T \DefinedEqual \sum_{k = 0}^{T - 1} \frac{d_k}{d_{k + 1}}
    = A_T + \sum_{k = 0}^{T - 1} \left(1 - \frac{d_k}{d_{k + 1}}\right) A_k
    \leq A_T^* (1 + T - S_T),
  \end{align*}
  where $A_T^* = \max_{0 \leq k \leq T} A_k \equiv \max_{1 \leq k \leq T} A_k$
  and we have used the fact that $(d_i)_{i = 0}^\infty$ is nondecreasing.
  Hence,
  \begin{align*}
    A_T^* \geq \frac{S_T}{1 + T - S_T}.
  \end{align*}
  Applying now the AM-GM inequality and
  denoting $\gamma_T = (\frac{d_0}{d_{T}})^{\frac{1}{T}}$
  ($\in \OpenClosedInterval{0}{1}$), we can further estimate
  $S_T \geq T \gamma_T$, giving us
  \begin{align*}
    A_T^* \geq \frac{T \gamma_T}{1 + T (1 - \gamma_T)}.
  \end{align*}
  Thus,
  \begin{align*}
    \min_{1 \le k \le T} \frac{d_k}{\sum_{i = 0}^{k - 1} d_i} =
    \frac{1}{A_T^*} & \leq
    \frac{\frac{1}{\gamma_T} \left(1 + T (1 - \gamma_T)\right)}{T}.
  \end{align*}
  Estimating further
  $T (1 - \gamma_T) \leq -T \log \gamma_T \equiv \log \frac{1}{\gamma_T^T}$
  and substituting the definition of $\gamma_T$, we get the claim.
\end{proof}

The following lemma is a classical result from
\cite[Lemma~3.2.1]{Nesterov2018LecturesOnCVX}.
\begin{lemma}\label{lem:omega-upper}
  Let $x \in \EffectiveDomain f$ be such that $f$ is subdifferentiable on $x$.
  Then, we have $f(x) - f^* \le \omega(v(x))$,
  where $\omega(\cdot)$ and $v(\cdot)$ are defined as in
  \eqref{eq:v-definition} and \eqref{eq:omega-definition}
  $($with $\nabla f(x)$ being an arbitrary subgradient from
  $\Subdifferential f(x)$$)$.
\end{lemma}

\begin{proof}
  Let $\bar{x}$ denote the orthogonal projection of $x^*$
  onto the supporting hyperplane
  \mbox{$\SetBuilder{y}{\inp{\nabla f(x)}{x - y} = 0}$}:
  \begin{align*}
    \bar{x} = x^* + v(x)\frac{B^{-1} \nabla f(x)}{\dnorm{\nabla f(x)}}.
  \end{align*}
  Then, $\inp{\nabla f(x)}{\bar{x} - x} = 0$, and $\norm{\bar{x} - x^*} = v(x)$.
  Therefore,
  \begin{align*}
    f(x)
    \leq f(\bar{x}) + \inp{\nabla f(x)}{\bar{x} - x} = f(\bar{x}),
  \end{align*}
  and hence,
  \begin{align*}
    f(x) - f^* & \leq
    f(\bar{x}) - f^* \le \omega(\norm{\bar{x} - x^*}) = \omega(v(x)).
    \qedhere
  \end{align*}
\end{proof}

\begin{lemma}\label{lem:bound-induction}
  Consider the nonnegative sequence $(d_k)_{k = 0}^\infty$
  that satisfies, for each $k \geq 0$,
  \begin{align*}
    d_{k + 1} \leq \max\{ d_k, R + \gamma d_k \},
  \end{align*}
  where $0 \leq \gamma < 1$ and $R \geq 0$ are certain constants.
  Then, for any $k \geq 0$, we have
  \begin{align*}
    d_k \leq \max\left\{\frac{1}{1 - \gamma} R, d_0\right\}.
  \end{align*}
\end{lemma}

\begin{proof}
  We use induction to prove that $d_k \leq D$
  for a certain constant $D$ to be determined later.
  To ensure that this relation holds for $k = 0$,
  we need to choose $D \geq d_0$. Let us now suppose that our
  relation has already been proved for some $k \geq 0$ and let us
  prove it for the next index $k + 1$.
  Using the induction hypothesis and the given inequality, we obtain
  \begin{align*}
    d_{k + 1}
    \leq
    \max\{ d_k, R + \gamma d_k \}
    \leq
    \max\{ D, R + \gamma D \}.
  \end{align*}
  To prove that the right-hand side is $\leq D$,
  we need to ensure that $R + \gamma D \leq D$,
  which means that we need to choose $D \geq \frac{1}{1 - \gamma} R$.
  Combining this requirement with that from the base of induction,
  we see that we can choose $D = \max\{\frac{1}{1 - \gamma} R, d_0\}$.
\end{proof}

\section{Proof of Theorem~\ref{thm:main-theorem}}
\label{sec:main-theorem-v}

\begin{lemma}\label{lem:da-convergence}
  In \cref{alg:da}, for any $1 \leq k \leq T$, it holds that
  \begin{align*}
    \sum_{i = 0}^{k - 1} a_i \inp{g_i}{x_i - x^*}
    + \frac{\beta_{k}}{2}\norm{x_k - x^*}^2 \le
    \frac{\beta_k}{2}\norm{x_0 - x^*}^2 + \sum_{i = 0}^{k -
      1}\frac{a_i^2}{2\beta_i}\dnorm{g_i}^2,
  \end{align*}
  where $\beta_0$ is an arbitrary coefficient in
  $\OpenClosedInterval{0}{\beta_1}$.
\end{lemma}

\begin{proof}
  For any $0 \leq k \leq T$, define the function $\psi_k(x)$ as follows:
  \begin{align*}
    \psi_k(x) \DefinedEqual \sum_{i = 0}^{k - 1} a_i \inp{g_i}{x - x_i} +
    \frac{\beta_k}{2} \norm{x - x_0}^2,
  \end{align*}
  so that $\psi_0(x) = \frac{\beta_0}{2} \norm{x - x_0}^2$ (with $\beta_0$ as
    defined in the statement). Note that
  $\psi_k$ is a $\beta_k$-strongly convex function and $x_k$ is its minimizer.
  Hence, for any $x \in Q$ and $0 \leq k \leq T$, we have
  \begin{align}\label{eq:phi-strong-convexity}
    \psi_{k}(x)
    \geq \psi_{k}^* + \frac{\beta_{k}}{2} \norm{x - x_{k}}^2,
  \end{align}
  where $\psi_{k}^* \DefinedEqual \psi_{k}(x_{k})$. Consequently,
  \begin{align*}
    \psi_{k + 1}^* & = \psi_{k + 1}(x_{k + 1}) = \psi_{k}(x_{k + 1})
    + a_{k}\inp{g_{k}}{x_{k + 1} - x_{k}}
    + \frac{\beta_{k + 1} - \beta_{k}}{2}\norm{x_{k + 1} - x_0}^2 \\
    & \geq \psi_{k}^* + \frac{\beta_{k}}{2} \norm{x_{k + 1} - x_{k}}^2
    + a_{k}\inp{g_{k}}{x_{k + 1} - x_{k}}
    + \frac{\beta_{k + 1} - \beta_{k}}{2}\norm{x_{k + 1} - x_0}^2 \\
    & \geq \psi_{k}^* + \frac{\beta_{k}}{2}\norm{x_{k + 1} - x_{k}}^2
    + a_{k} \inp{g_{k}}{x_{k + 1} - x_{k}} \geq \psi_{k}^* -
    \frac{a_{k}^2}{2\beta_{k}}\dnorm{g_{k}}^2.
  \end{align*}
  Telescoping these inequalities and
  using the fact that $\psi_0^* = 0$, we obtain
  \begin{align*}
    \psi_{k}^* & \ge - \sum_{i = 0}^{k - 1} \frac{a_i^2}{2\beta_i}
    \dnorm{g_i}^2.
  \end{align*}
  Combining this inequality with the definition of $\psi_k$ and
  \cref{eq:phi-strong-convexity}, we thus obtain
  \begin{align*}
    \sum_{i = 0}^{k - 1} a_i \inp{g_i}{x^* - x_i}
    + \frac{\beta_k}{2} \norm{x_0 - x^*}^2 =
    \psi_{k}(x^*) & \ge \psi_{k}^* + \frac{\beta_{k}}{2} \norm{x_{k} - x^*}^2
    \\
    & \ge - \sum_{i = 0}^{k - 1} \frac{a_i^2}{2\beta_i} \dnorm{g_i}^2 +
    \frac{\beta_{k}}{2} \norm{x_{k} - x^*}^2.
  \end{align*}
  Rearranging, we get the
  claim.
\end{proof}

\begin{lemma}\label{lem:coefficient-placing}
  Consider \cref{alg:da} using the coefficients
  defined in \cref{eq:dada-stepsizes}.
  Then, the following inequality holds for all $1 \leq k \leq T$:
  \begin{align*}
    \sum_{i = 0}^{k - 1} \bar{r}_i v_i
    + \frac{c \sqrt{k + 1}}{2} D_k^2
    \leq \frac{c \sqrt{k + 1}}{2} D_0^2
    + \frac{\sqrt{k}}{c} \bar{r}_{k - 1}^2,
  \end{align*}
  where $D_k = \norm{x_k - x^*}$ and $v_i \DefinedEqual v(x_i)$.
\end{lemma}

\begin{proof}
  Applying \cref{lem:da-convergence} and the definition of $v_i$, we obtain
  \begin{align*}
    \sum_{i = 0}^{k - 1} a_i v_i \dnorm{g_i} + \frac{\beta_k}{2} D_k^2
    \leq \frac{\beta_k}{2} D_0^2 + \sum_{i = 0}^{k - 1} \frac{a_i^2}{2\beta_i}
    \dnorm{g_i}^2.
  \end{align*}
  Substituting our choice of the coefficients
  given by \eqref{eq:dada-stepsizes}, we get
  \begin{align*}
    \sum_{i = 0}^{k - 1} \bar{r}_i v_i
    + \frac{c \sqrt{k + 1}}{2} D_k^2
    & \leq \frac{c \sqrt{k + 1}}{2} D_0^2
    + \frac{1}{2c}
    \sum_{i = 0}^{k - 1} \frac{\bar{r}_i^2}{\sqrt{i + 1}}
    \leq \frac{c \sqrt{k + 1}}{2} D_0^2
    + \frac{\sqrt{k}}{c} \bar{r}_{k - 1}^2,
  \end{align*}
  where we have used the fact that $\bar{r}_k$ is nondecreasing and
  $\sum_{i = 0}^{k - 1} \frac{1}{\sqrt{i + 1}} \le 2 \sqrt{k}$.
\end{proof}

\begin{lemma}\label{lem:r-upper-d}
  Consider \cref{alg:da} using the coefficients
  defined in \cref{eq:dada-stepsizes} and assume that $c > \sqrt{2}$.
  Then, we have the following inequalities for all $0 \leq k \leq T$:
  \[
    \bar{r}_k \leq \bar{D},
    \qquad
    D_k \leq D_0 + \frac{\sqrt{2}}{c} \bar{D},
  \]
  where
  $
    \bar{D}
    \DefinedEqual
    \max \bigl\{ \bar{r}, \frac{2 c}{c - \sqrt{2}} D_0 \bigr\}
  $
  and $D_k \DefinedEqual \Norm{x_k - x^*}$.
\end{lemma}

\begin{proof}
  Both bounds are clearly valid for $k = 0$, so it suffices to consider only
  the case when $1 \leq k \leq T$.

  Applying \cref{lem:coefficient-placing}, dropping the
  nonnegative $\bar{r}_i v_i$ from the left-hand side and rearranging,
  we obtain
  \[
    D_k^2
    \leq D_0^2 + \frac{2 \sqrt{k}}{c^2 \sqrt{k + 1}} \Bar{r}_{k - 1}^2
    \leq D_0^2 + \frac{2}{c^2} \Bar{r}_{k - 1}^2.
  \]
  Consequently,
  \begin{equation}\label{eq:dk-inequality}
    D_k \le D_0 + \frac{\sqrt{2}}{c} \Bar{r}_{k - 1}.
  \end{equation}
  Therefore,
  \[
    r_k \equiv \norm{x_k - x_0}
    \leq
    D_k + D_0 \leq 2 D_0 + \frac{\sqrt{2}}{c} \Bar{r}_{k - 1}.
  \]
  Hence,
  \begin{align*}
    \bar{r}_k
    \equiv
    \max\{ \bar{r}_{k - 1}, r_k \}
    \leq
    \max\biggl\{ \bar{r}_{k - 1}, 2 D_0 + \frac{\sqrt{2}}{c} \bar{r}_{k - 1}
    \biggr\}.
  \end{align*}
  Since $k \geq 1$ was allowed to be arbitrary, we can apply
  \cref{lem:bound-induction} to conclude that
  \begin{align*}
    \bar{r}_k & \leq
    \max \biggl\{\bar{r},
      \frac{2}{1 - \frac{\sqrt{2}}{c}} D_0
    \biggr\} = \max \biggl\{\bar{r},
      \frac{2 c}{c - \sqrt{2}} D_0
    \biggr\}
    \equiv
    \bar{D}.
  \end{align*}
  This proves the first part of the claim.

  Substituting the already proved bound on~$\bar{r}_k$ into
  \cref{eq:dk-inequality}, we obtain the claimed upper bound on~$D_k$.
\end{proof}
We are now ready to prove the main result.
\begin{proof}[Proof of \cref{thm:main-theorem}]\label{prf:main-theorem}
  Let $T \geq 1$ be arbitrary.
  According to \cref{lem:omega-upper} and the fact that $\omega(\cdot)$ is
  nondecreasing, we can write
  \begin{align*}
    f(x_T^*) - f^* = \min_{0 \le k \le T - 1} [f(x_k) - f^*] \leq
    \min_{0 \le k \le T - 1} \omega(v_k) = \omega(v_T^*),
  \end{align*}
  where $v_k \DefinedEqual v(x_k)$ and
  $v_T^* \DefinedEqual \min_{0 \leq k \leq T - 1} v_k$.
  This proves the first part of the claim.

  Let us now estimate the rate of convergence of~$v_T^*$.
  To that end, let us fix an arbitrary $1 \leq k \leq T$.
  In view of \cref{lem:coefficient-placing}, we have
  \begin{align*}
    \sum_{i = 0}^{k - 1} \bar{r}_i v_i & \leq
    \frac{c \sqrt{k + 1}}{2} (D_0^2 - D_k^2) +
    \frac{\sqrt{k}}{c} \bar{r}_{k - 1}^2,
  \end{align*}
  where $D_k = \norm{x_k - x^*}$.
  Note that
  \begin{align*}
    D_0^2 - D_k^2
    \equiv \norm{x_0 - x^*}^2 - \norm{x_k - x^*}^2
    &= (\norm{x_0 - x^*} - \norm{x_k - x^*})
    (\norm{x_0 - x^*} + \norm{x_k - x^*}) \\
    & \leq 2 \norm{x_k - x_0}\norm{x_0 - x^*} \equiv 2 r_k D_0.
  \end{align*}
  Therefore, we can continue as follows:
  \begin{align*}
    \sum_{i = 0}^{k - 1} \Bar{r}_i v_i
    &\leq
    c \sqrt{k + 1} r_k D_0
    +
    \frac{\sqrt{k}}{c} \bar{r}_{k - 1}^2
    \leq
    \Bigl( c D_0 + \frac{1}{c} \bar{r}_{k - 1} \Bigr) \sqrt{k + 1} \, \bar{r}_k
    \\
    &\leq
    \Bigl( c D_0 + \frac{1}{c} \bar{D} \Bigr) \sqrt{k + 1} \, \bar{r}_k
    =
    D \sqrt{ \frac{k + 1}{2} } \bar{r}_k,
  \end{align*}
  where the second inequality is due to the fact that
  $\bar{r}_k = \max\{ \bar{r}_{k - 1}, r_k \}$,
  the final inequality is due to \cref{lem:r-upper-d},
  and the constants $\bar{D}$ and~$D$ are as defined in the statement.
  Hence,
  \begin{align*}
    v_k^* \equiv \min_{0 \leq i \leq k - 1} v_i \leq
    \frac{\sum_{i = 0}^{k - 1}\Bar{r}_i v_i}{\sum_{i = 0}^{k - 1} \Bar{r}_i}
    \leq
    \frac{\Bar{r}_{k}}{\sum_{i = 0}^{k - 1} \Bar{r}_i}
    D \sqrt{\frac{k + 1}{2}}.
  \end{align*}
  Letting now
  $
    k^* = \argmin_{1 \le k \le T}
    \frac{\Bar{r}_k}{\sum_{i = 0}^{k - 1} \Bar{r}_i}
  $
  and using \cref{lem:seq-min-sum}, we obtain
  \[
    v_T^* \leq v_{k^*}^*
    \leq \frac{D \sqrt{\frac{k^* + 1}{2}}}{T}
    \left(\frac{\Bar{r}_T}{\Bar{r}}\right)^{\frac{1}{T}} \log \frac{e
      \Bar{r}_T}{\Bar{r}}
    \leq \frac{D}{\sqrt{T}}
    \left(\frac{\bar{D}}{\Bar{r}}\right)^{\frac{1}{T}}
    \log \frac{e \bar{D}}{\Bar{r}},
  \]
  where we have used the fact that $k^* + 1 \leq T + 1 \leq 2 T$
  (since $1 \leq k^* \leq T$) and that $\bar{r}_T \leq \bar{D}$
  (see \cref{lem:r-upper-d}). This proves
  \cref{eq:main-theorem:convergence-rate} in the case when
  $T \geq \log \frac{\bar{D}}{\Bar{r}}$ since then we can further bound
  $
    (\frac{\bar{D}}{\bar{r}})^{\frac{1}{T}}
    \equiv
    \exp(\frac{1}{T} \log \frac{\bar{D}}{\Bar{r}})
    \leq
    e
  $.

  On the other hand, by the definition of~$v_k$ and \cref{lem:r-upper-d}, we
  always have the following trivial inequality for any $0 \leq k \leq T - 1$:
  \[
    v_k
    \equiv \frac{\inp{\nabla f(x_k)}{x_k - x^*}}{\dnorm{\nabla f(x_k)}}
    \leq
    D_k
    \leq
    D_0 + \frac{\sqrt{2}}{c} \bar{D}
    \leq
    D.
  \]
  This means that \cref{eq:main-theorem:convergence-rate} is also satisfied
  in the case when $T \leq \log \frac{\bar{D}}{\bar{r}}$ since then
  $
    \frac{e D}{\sqrt{T}} \log \frac{e \bar{D}}{\bar{r}}
    \geq
    \frac{D}{\sqrt{T}} \log \frac{\bar{D}}{\bar{r}}
    \geq
    D \sqrt{T}
    \geq
    D
  $
  (we still consider $T \geq 1$). The proof of
  \cref{eq:main-theorem:convergence-rate} is now finished.

  The final part of the claim readily follows from
  \cref{eq:main-theorem:convergence-rate}.
\end{proof}

\section{How to choose the constant $c$}\label{sec:discussion-on-c}
According to \cref{thm:main-theorem}, our method converges for any
$c > \sqrt{2}$. However, the choice of $c$ can influence the constant factor
in the complexity of DADA. Hence, our goal here is to find the optimal constant
$c$ that minimizes $T_v(\delta)$.
To determine this $c$, let $\bar{r}$ be sufficiently small, so that
\[
  \bar{D}
  \equiv
  \max\Bigl\{ \bar{r}, \frac{2c}{c - \sqrt{2}}D_0 \Bigr\}
  =
  \frac{2c}{c - \sqrt{2}}D_0.
\]
Then, disregarding the logarithmic factors,
due to their minimal impact on the complexity of our method,
we can determine the optimal constant $c$ that minimizes
\begin{align*}
  D
  \equiv
  \sqrt{2} \Bigl( c D_0 + \frac{1}{c} \bar{D} \Bigr)
  =
  \sqrt{2} \Bigl( c + \frac{2}{c - \sqrt{2}} \Bigr) D_0.
\end{align*}
This is the value
\begin{align}\label{eq:optimal-c}
  c = 2\sqrt{2}.
\end{align}
For this optimal choice of $c$, we get $\bar{D} = \max \{\bar{r}, 4 D_0\}$ and
$D = 4 D_0 + \frac{1}{2} \bar{D}$, so the complexity of our method given by
\cref{thm:main-theorem} is
\[
  T_v(\delta) =
  \frac{e^2 (4 D_0 + \frac{1}{2} \bar{D})^2}{\delta^2}
  \log^2 \frac{e \bar{D}}{\Bar{r}}.
\]

\section{Convergence of DADA on Various Problem Classes}
\label{sec:convergence-proof-examples}

In this section, we analyze the complexity of DADA across
different problem classes. To achieve this, we first establish
bounds on the growth function:
\begin{align*}
  \omega(t) = \max_{x \in \Ball{x^*}{t}} f(x) - f^*,
\end{align*}
and determine the threshold $t$ such that $\omega(t) \leq \epsilon$
for a given $\epsilon$. Subsequently, we combine these results with the
complexity bound $T(\delta)$ derived in \cref{thm:main-theorem},
enabling us to estimate the oracle complexity of DADA for finding an
$\epsilon$-solution in terms of the function residual.

\section{Additional Experiments}\label{sec:additional-experiments}

\paragraph{Softmax function.}
Our first test problem is
\begin{align}\label{eq:softmax}
  \min_{x \in \R^d} \Biggl\{f(x) \DefinedEqual
    \mu \log\left(
      \sum_{i = 1}^{n} \exp\left[\dfrac{\InnerProduct{a_{i}}{x} - b_i}{\mu}\right]
    \right)\Biggr\},
\end{align}
where $a_i \in \R^d$, and $b_i \in \R$ for all $1 \leq i \leq n$, and $\mu >
  0$.
This function can be viewed as a smooth approximation of
$\max_{1 \leq i \leq n} [\inp{a_i}{x} - b_i]$
\cite{nesterov2005smooth}.

To generate the data for our problem, we proceed as follows.
First, we generate i.i.d.\ vectors $\hat{a}_i$ with components
uniformly distributed in the interval $[-1, 1]$ for $i = 1, \dots, n$,
and similarly for the scalar values $b_i$.
Using this data, we form the preliminary version of our function, $\hat{f}$.
We then compute $a_i = \hat{a}_i - \nabla \hat{f}(0)$ and use the obtained
$(a_i, b_i)$ to define our
function $f$. This way of generating the data ensures that $x^* = 0$ is a
solution of our problem.

\begin{figure}[tb]
  \centering
  \includegraphics[width=\textwidth]{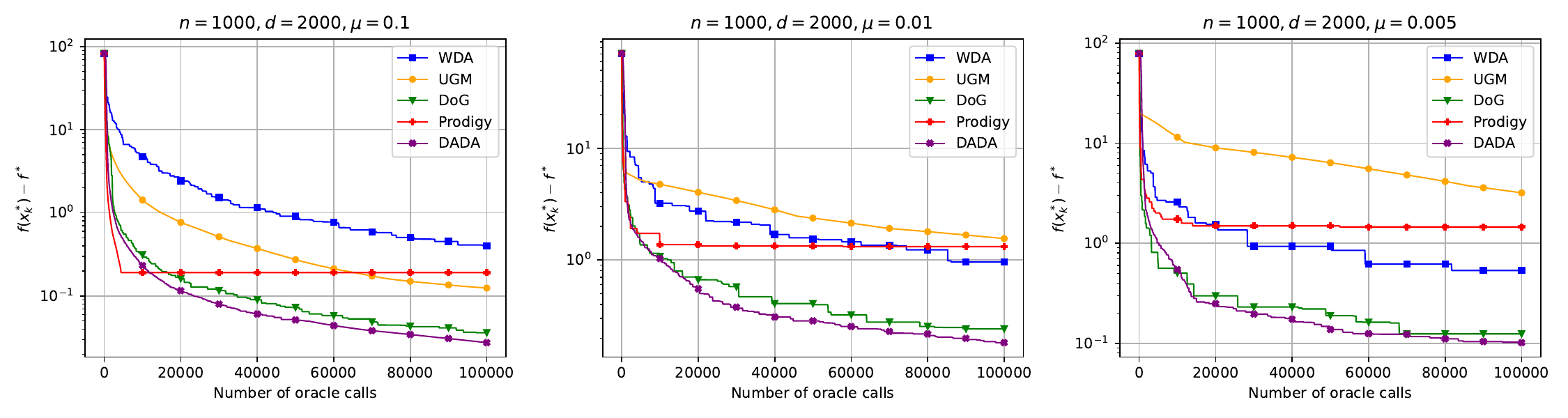}
  \caption{
    Comparison of different methods on the Softmax function.
  }
  \label{fig:softmax-residual}
\end{figure}
The results are shown in \cref{fig:softmax-residual},
where we fix $n = 10^3$ and $d = 2 n$, and consider
different values of $\mu \in \Set{0.1, 0.01, 0.005}$.
As we can see, most methods exhibit similar
performance for $\mu = 0.1$ except for Prodigy
which stops converging after a few initial iterations. This issue,
along with a decline in performance for UGM, persists as $\mu$ decreases,
whereas DADA, DoG, and WDA remain largely unaffected.
Notably, DoG performs very similarly to DADA, which we
hypothesize is primarily due to the similarity in estimating~$D_0$.
\begin{figure}[tb]
  \centering
  \includegraphics[width=\textwidth]{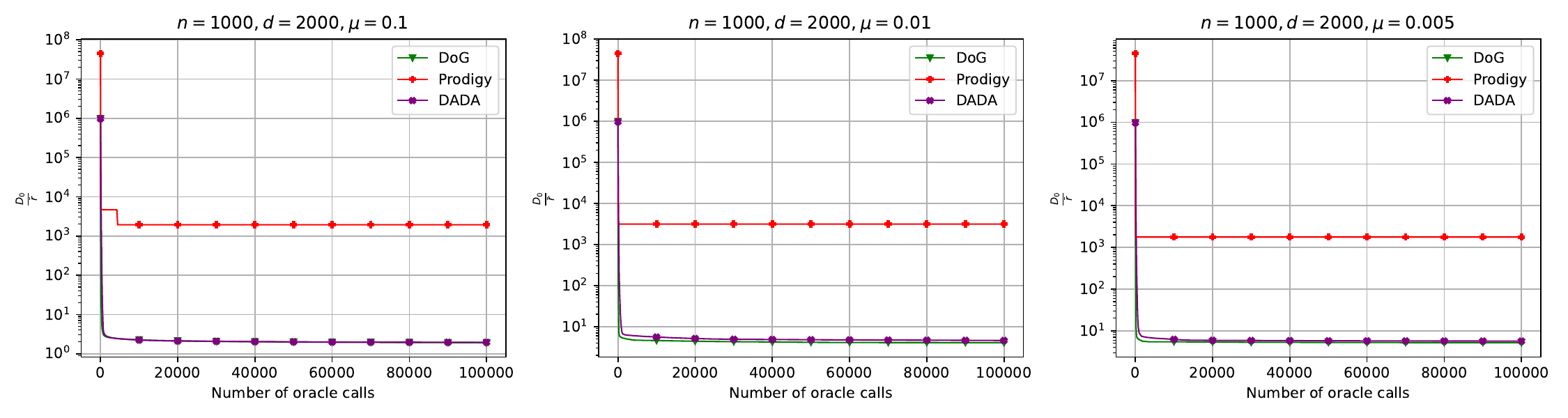}
  \caption{
    The ratio $\frac{D}{\bar{r}_t}$
    for the Softmax function with different optimal points $x^*$.
  }
  \label{fig:softmax-estimate-error}
\end{figure}

Additionally, \cref{fig:softmax-estimate-error} illustrates the
ratio between $D_0$ and $\bar{r}$,
showing the estimation error of Prodigy, DoG, and DADA
throughout the optimization process.
For Prodigy, we use \(\frac{D_0}{d_{\max}}\) to generate the plot.
The figure demonstrates that DADA and DoG exhibit similar
behavior in estimating $D_0$, despite employing different
update methods—Dual Averaging and Gradient Descent, respectively.
However, Prodigy appears to encounter challenges in estimating~$D_0$
as its estimation stabilizes at a relatively large value.

\paragraph{H\"older-smooth function.}
Let us consider the following \emph{polyhedron feasibility problem}:
\begin{align}\label{eq:holder-smooth-sample}
  f^* \DefinedEqual \min_{x \in \R^d} \Bigl\{
    f(x) \DefinedEqual \frac{1}{n} \sum_{i = 1}^{n}
    \PositivePart{\InnerProduct{a_i}{x} - b_i}^q
  \Bigr\},
\end{align}
where $a_i, b_i \in \R^d$, $q \in \ClosedClosedInterval{1}{2}$,
and $\PositivePart{\tau} = \max(0, \tau)$.
This problem can be interpreted as finding a point $x^* \in \R^d$ lying
inside the polyhedron
$P = \SetBuilder{x}{\inp{a_i}{x} \leq b_i, \ i = 1, \ldots, n}$.
Such a point exists if and only if $f^* = 0$.

Observe that $f$ in problem~\eqref{eq:holder-smooth-sample} is H\"older-smooth
with parameter $\nu = q - 1$. Therefore, by varying
$q \in \ClosedClosedInterval{1}{2}$, we can check the robustness of different
methods to the smoothness level of the objective function.

The data for our problem is generated randomly,
following the procedure in \cite{Rodomanov2024UniversalityAdaGrad}.
First, we sample $x^*$ uniformly from the sphere of radius $0.95 R$
centered at the origin.
Next, we generate i.i.d.\ vectors $a_i$ with components
uniformly distributed in $\ClosedClosedInterval{-1}{1}$.
To ensure that $\InnerProduct{a_n}{x^*} < 0$,
we invert the sign of $a_n$ if necessary.
We then sample positive reals $s_i$ uniformly from
$\ClosedClosedInterval{0}{-0.1 c_{\min}}$, where
$c_{\min} \DefinedEqual \min_i \InnerProduct{a_i}{x^*} < 0$,
and set $b_i = \InnerProduct{a_i}{x^*} + s_i$.
By construction, $x^*$ is a solution to the problem with $f^* = 0$.

\begin{figure}[tb]
  \centering
  \includegraphics[width=\textwidth]{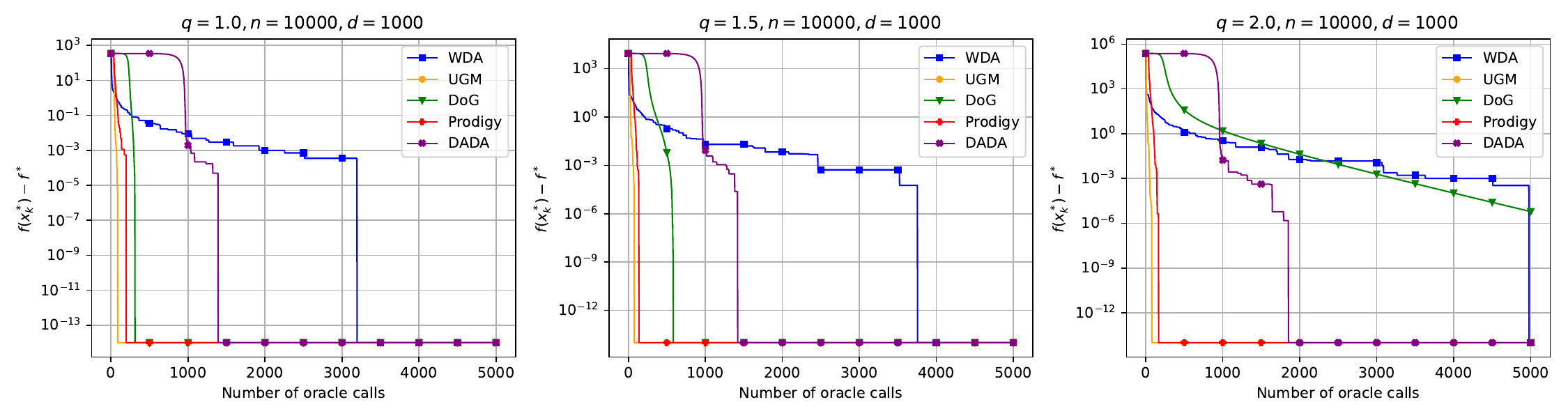}
  \caption{
    Comparison of different methods on the polyhedron feasibility problem.
  }
  \label{fig:holder-smooth-residual}
\end{figure}

We select $n = 10^4$, $d = 10^3$, $R = 10^3$ and consider different values of
$q \in \Set{1, 1.5, 2}$. As shown in \cref{fig:holder-smooth-residual},
as $q$ increases and approaches 2,
the performance of DoG significantly declines.
However, DADA, Prodigy, and UGM demonstrate similar performance
regardless of the choice of $q$.

\paragraph{Comparison of different initial estimates of the distance.}
In this experiment, we evaluate the sensitivity of DADA to the choice
of the initial point $x_0$.
We consider the same Softmax function as in \eqref{eq:softmax}
with $n = 10^3$, $d = 2 n$, and $\mu \in \Set{0.5, 0.1, 0.01}$.

\begin{figure}[tb]
  \centering
  \begin{subfigure}[b]{0.32\textwidth}
    \centering
    \includegraphics[width=\textwidth]{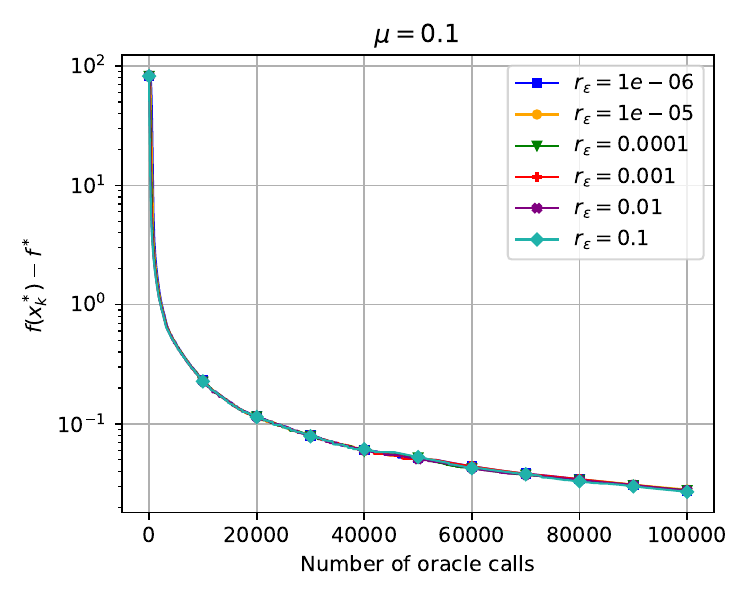}
    \label{fig:plot1}
  \end{subfigure}
  \begin{subfigure}[b]{0.32\textwidth}
    \centering
    \includegraphics[width=\textwidth]{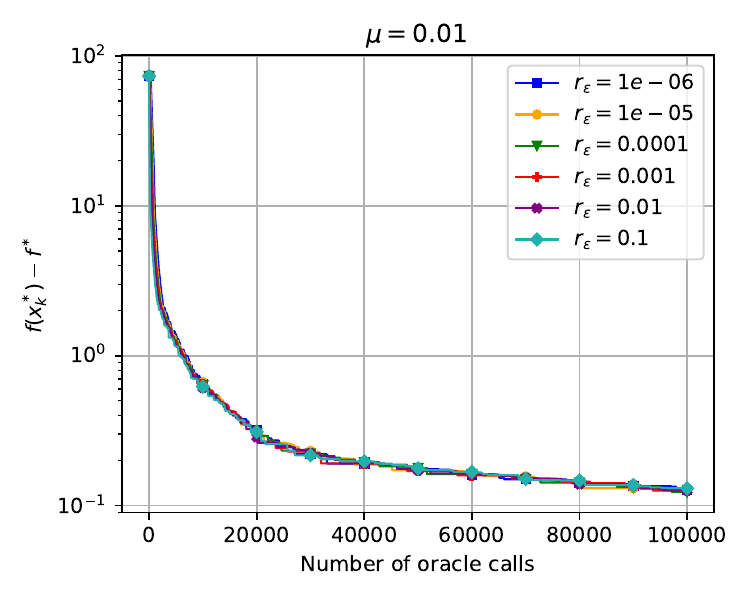}
    \label{fig:plot2}
  \end{subfigure}
  \begin{subfigure}[b]{0.32\textwidth}
    \centering
    \includegraphics[width=\textwidth]{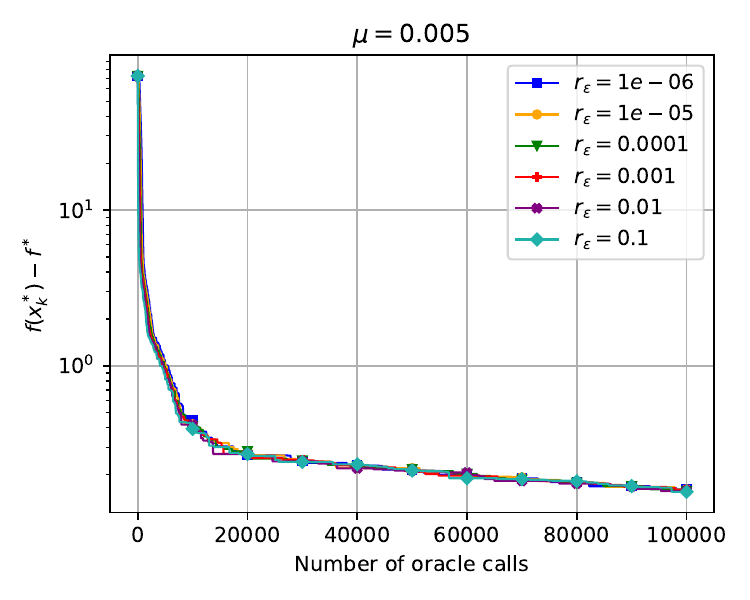}
    \label{fig:plot3}
  \end{subfigure}

  \caption{Comparison of different initial estimates of the distance on the
    Softmax function with
    different values of $\mu$.}
  \label{fig:different-initial-point}
\end{figure}

The results are shown in \cref{fig:different-initial-point},
where we consider $\delta \in \Set{10^{-1}, \ldots, 10^{-6}}$.
As we can see, the choice of $\delta$ does not affect the performance of
DADA, which consistently achieves similar performance across all tested values.

\end{document}